\documentclass{amsart}[11pt]
\usepackage{myAMSpackages}
\usepackage{Universal}
\usepackage{MyAMSEnvironments}
\usepackage{graphicx,subfigure,epsfig}

\newcommand{\tor}{\text{Tor}}
\newcommand{\Tor}{\text{TOR$_2$}}
\newcommand{\torr}[1][r]{\text{Tor}^{(#1)}}
\newcommand{\tr}[1][0]{\ensuremath{\text{tr}_{(#1)}}\xspace}

\newcommand{\Ric}{\ensuremath{\text{Rc}^s}\xspace}
\newcommand{\srm}{\ensuremath{\text{Rm}^s}\xspace}
\newcommand{\srt}{\ensuremath{\C\Tor}\xspace}

\newcommand{\V}[1][*]{\ensuremath{V^{(#1)}}\xspace}
\newcommand{\hV}[1][*]{\ensuremath{\widehat{V}^{(#1)}}\xspace}
\newcommand{\nb}[1][r]{\nabla^{(#1)}}

\newcommand{\bn}{\overline{\nabla}}
\newcommand{\bR}{\overline{\text{Rm}}}

\newcommand{\B}[1][0]{\ensuremath{B^{(#1)}}\xspace}

\newcommand{\lp}[1][0]{\triangle_{(#1)}}

\newcommand{\C}{\mathscr{C}}

\begin{document}
\title{Connections and Curvature in subRiemannian geometry}
\author{Robert K. Hladky}
\address{North Dakota State University Dept. \#2750, 
PO Box 6050, 
Fargo ND 58108-6050}
\email{robert.hladky@ndsu.edu}
\urladdr{http://www.ndsu.edu/pubweb/~hladky/}
\subjclass{53C17, 53C05}
\keywords{SubRiemannian geometry, Carnot-Carath\'eodory geometry, curvature,  Bochner Formula, Bianchi Identity, Killing fields}

\begin{abstract}
For a subRiemannian manifold and a given Riemannian extension of the metric, we define a canonical global  connection. This connection coincides with both the Levi-Civita connection on Riemannian manifolds and the Tanaka-Webster connection on strictly pseudoconvex CR manifolds.  We define a notion of normality generalizing Tanaka's notion for CR manifolds to the subRiemannian case. Under the assumption of normality, we construct local frames that simplify computations in a manner analogous to Riemannian normal coordinates. We then use these frames to establish Bianchi Identities and symmetries for the associated curvatures. Next we study subRiemannian notions of the Ricci curvature and horizontal Laplacian, establishing general Bochner type identities. Finally we explore subRiemannian generalizations of the Bonnet-Myers theorem, providing some new results and some new proofs and interpretations of existing results in the literature.
\end{abstract}

\maketitle

\section{Introduction}\setS{ID}

A fundamental tool in Riemannian geometry is the Levi-Civita connection. As the device which permits us to glue local differential equations into global ones, it is the key ingredient in most modern descriptions of curvature and geodesics and underlies many computational methods in differential geometry. The Tanaka-Webster connection, (\cite{Tanaka}, \cite{Webster}) plays a similar role in the study of strictly pseudoconvex CR manifolds. 

There has been much recent effort to define such geometrically useful connections in sub-Riemannian geometry. All work, including this one, has operated under the assumption that the subRiemannian metric on the horizontal bundle has  been extended to a Riemannian metric on the whole space. This allows us to define a vertical bundle. Previous work has been inherently local, depending on a choice of frame.  Usually some additional geometric  and topological restrictions have been required. In \cite{HP2}, \cite{HP4} a subRiemannian connection was defined under the assumption of a global frame for the vertical bundle. In \cite{BaudoinGarofalo}, a connection is defined under a strong tensorial condition, referred to as strict normality in this paper and the assumption of the existence of a  frame of vertical Killing fields.  All of these examples required a \textit{a priori} choice of frame for the vertical bundle and so do not define global connections in general.

This lack of a global covariant derivative scheme means that the study of the relationships between subelliptic PDE and subRiemannian manifolds has been by necessity local in nature.  Recently there has been some effort addressing this need. In particular in \cite{BaudoinGarofalo}, several global curvature results such as Myer's theorem have been extended to certain step 2 subRiemannian manifolds.

In this paper, we propose a new globally defined connection to facilitate this process. We shall work assumption that a global complement to the horizontal bundle has been chosen. For any Riemannian metric extending the subRiemannian metric and preserving this decomposition, we shall define a canonical, global metric compatible connection such that the horizontal and vertical bundles are parallel.  In the special cases of Riemannian and strictly pseudoconvex pseudohermitian manifolds, this connection will coincides with the Levi-Civita and the Tanaka-Webster connections respectively. Furthermore any covariant derivative of any horizontal vector field will be seen to be independent of the choice of Riemannian extension. Thus for a subRiemannian manifold with vertical complement, there is a canonical method for taking covariant derivatives of horizontal vector fields.

In section 2 we define the connection and explore its basic properties and how they relate to bracket structures of the underlying horizontal and vertical bundles. We introduce a tool similar to Riemannian normal coordinates, to aid computation. In section 3 we consider the associated curvature tensors and their symmetries. SubRiemannian equivalents of the Bianchi identities are introduced and proved. In section 4 we establish some Bochner-type formulas for general subRiemannian manifolds and show how the analytic framework developed by Baudoin and Garofalo generalizes to the category strictly normal subRiemannian manifolds. In section 5, we compare the subRiemannian connection to the Levi-Civita connection for metric extensions. We then use this to provide a new interpretation and proof of an existing subRiemannian Bonnet-Myers theorem as well as providing new results in this vein.

\section{SubRiemannian manifolds and the Canonical Connection}\setS{CC}

There are several subtly different notions of subRiemannian manifolds in the literature. In this paper, we shall use the following definition:
\bgD{SR}
A subRiemannian manifold is a smooth manifold $M$, a smooth constant rank distribution $HM \subset TM$ and a smooth inner product $\aip{\cdot}{\cdot}{}$ on $HM$. The bundle $HM$ is known as the horizontal bundle.
\enD

We should remark here that we are not assuming any conditions on the horizontal bundle other than constant rank. In particular in this paper, unless otherwise stated, we are not even assuming that it bracket generates.

\bgD{SRC}
A subRiemannian manifold with complement, henceforth sRC-manifold, is a subRiemannian manifold together with a smooth bundle $VM$ such that $HM \oplus VM = TM$. The bundle $VM$ is known as the vertical bundle.

Two sRC-manifolds  $M,N$are sRC-isometric if there exists a diffeomorphism $\pi\colon M \to N$ such that $\pi_*  HM =HN$, $\pi_* VM = VN$ and $\aip{\pi_* X}{\pi_* Y}{N} = \aip{X}{Y}{M}$ for all horizontal vectors  $X,Y$.
\enD

\bgD{REr} 
A sRC-manifold $(M,HM,VM,\aip{\cdot}{\cdot}{})$ is  $r$-graded if there are smooth constant rank bundles $\V[j]$, $0<j \leq r $, such that
\[ VM = \V[1] \oplus \dots \oplus \V[r] \] \bgE{grading} HM \oplus \V[j] \oplus \left[HM, \V[j]\right]  \subseteq HM \oplus \V[j] \oplus \V[j+1]  \enE
for all $0\leq j  \leq r$.  Here we have adopted the convention that $\V[0]=HM$ and $\V[k]=0$ for $k>r$.

The grading is $j$-regular if  
\bgE{regular}
HM \oplus \V[j] \oplus \left[HM, \V[j]\right]  =  HM \oplus \V[j] \oplus \V[j+1] 
\enE  and equiregular if  is $j$-regular for all $0 \leq j \leq r$. 

A metric extension for an $r$-graded vertical complement is a Riemannian metric $g$ of $\aip{\cdot}{\cdot}{}$ that makes the split
\[ TM =HM   \bigoplus\limits_{1\leq j \leq r}  \V[j] \]
orthogonal.

\enD

For convenience of notation, we shall denote a section of $\V[k]$ by $X^{(k)}$ and set
\[ \hV[j] = \bigoplus\limits_{k \ne j} \V[k]. \]
If a metric extension has been chosen then $\hV[j]= \left( \V[j] \right)^\perp$.
the orthogonal complement of $\V[j]$. For convenience, we shall often also extend the notation $\aip{\cdot}{\cdot}{}$ to whole tangent space using it interchangeably with $g$.

\bgR{subgrading}
Every sRC-manifold that admits an $r$-grading also admits $k$-gradings for all $1 \leq k <r$ by setting
\[ \tilde{V}^{(j)} = \V[j] \quad  0\leq j <k, \qquad  \tilde{V}^{(k)} = \bigoplus\limits_{j \geq k} \V[j] \]
\enR

\bgD{Basic}
The unique $1$-grading on each sRC-manifold,
\[ \V[1] =VM \]
is known as the basic grading.
\enD

\bgX{Carnot}
A Carnot group (of step $r$)  is a Lie group, whose Lie algebra $\mathfrak{g}$ is stratified in the sense that
\[  \mathfrak{g} = \mathfrak{g}_0 \oplus \dots \mathfrak{g}_{r-1}, \qquad [ \mathfrak{g}_0,\mathfrak{g}_j]=\mathfrak{g}_{j+1} \quad  j=1\dots r, \quad \mathfrak{g}_{r}=0\]
together with a left-invariant metric $\aip{\cdot}{\cdot}{}$ on $HM$, the left-translates of $\mathfrak{g}_0$.

The vertical bundle $VM$ consists of the left-translates of $ \mathfrak{g}_1 \oplus \dots \mathfrak{g}_{r-1}$.  In addition to the basic grading, there is then a natural equiregular $r-1$-grading defined by setting $\V[j]$ to be the left-translates of $\mathfrak{g}_j$.
\enX

\bgD{B}If a metric extension $g$ has been chosen, we define
\[ B(X,Y,Z) = (\mathcal{L}_Z g)(X,Y) = Z g(X,Y) + g([X,Z],Y) + g([Y,Z],X)\]
for vector fields $X,Y,Z$
\enD
Unfortunately $B$ is not tensorial in general and so cannot be viewed as a map on vectors rather than vector fields. However,  we can  define a symmetric tensor $\B[j]$ by setting
\[  \B[j](X,Y,T) = B(X,Y,T)\]
 for  $X,Y \in \V[j]$, $T \in \hV[j]$ and declaring $\B[j]$ to be zero on the orthogonal complement of $\V[j] \times \V[j] \times \hV[j]$.  We can then contract these to tensors $C^{(j)} \colon TM \times TM \to \V[j]$ defined by
\bgE{C}
g(C^{(j)}(X,Y),Z^{(j)}) = \B[j](X,Z^{(j)},Y) 
\enE
Additionally, we can define $j$-traces, by
\begin{align*}
\tr[j] \B[j](Z) &= \sum \B[j](E^{(j)}_i,E^{(j)}_i,Z) \\
\end{align*}
where $\{E^{(j)}_i\}$ are (local) orthonormal frames  for $\V[j]$.

\bgD{Normal}
Suppose that $M$ is an $r$-graded sRC-manifold with metric extension $g$.
\begin{itemize}
\item The metric extension is $j$-normal with respect to the grading if $\B[j] \equiv 0$.
\item The metric extension is strictly normal with respect to the grading if it is $j$-normal for all $0\leq j \leq r$. 
\end{itemize} 
\enD

\bgR{indep}
The tensors $\B$ and $C^{(0)}$ depend only on the underlying sRC-structure and are independent of the choice of grading and metric extension. Thus the notion of $0$-normal is also independent of grading and metric. We shall say $VM$ is normal, if every metric extension and grading is $0$-normal.
\enR

\bgX{C3}
Let $M$ be the 4 dimensional Carnot group with Lie algebra induced by the global left invariant vector fields $X,Y,T,S$ with bracket structures
\[ [X,Y]=T, \qquad [X,T] =S \]
and all others being zero. Then $B(T,S,X) = -1$ with all others vanishing.  Now $M$ admits an equiregular $2$-grading defined by 
\[ \V[1] = \langle T \rangle, \qquad \V[2] = \langle S \rangle .\]
Let $g$ be the metric making the global frame orthonormal. Then $g$ is strictly normal with respect to this $2$-grading.

It should be remarked that this metric is not $1$-regular with respect to the basic grading.  For then we  get $\tilde{B}^{(0)} \equiv 0$ but $\tilde{B}^{(1)}(T,S,X)=-1$.  Thus the metric is $0$-normal but not strictly normal with respect to the basic grading.
\enX

\bgX{CarnotN}
Any step $r$ Carnot group with a bi-invariant metric extension is strictly normal with respect to the equiregular  $r-1$-grading, but is only $0$-normal with respect to the basic grading.
\enX

\bgX{spsic}
Let $(M,J,\eta)$ be a strictly pseudoconvex pseudohermitian manifold, (see \cite{Tanaka}) with characteristic vector field $T$ such that $\eta(T)=1$, $T \lrcorner d\eta =0$. The horizontal bundle $HM$ is defined to be the kernel of the 1-form $\eta$. An immediate consequence of the defining properties of $T$ is that $[T,HM] \subset HM$. When $J$ is extended to $TM$ by defining $JT=0$, the Levi metric 
\[ g(A,B) = d\eta(A,JB) + \eta(A) \eta(B)\]
can be viewed as an extension of the subRiemannian metric $\aip{X}{Y}{} = d\theta(X,JY)$  with  $VM =\langle T \rangle$.  As $VM$ is one dimensional, the basic grading is the only grading admitted and  since $[T,HM] \subset HM$ we see $\B[1]=0$ trivially. Thus the Levi metric is always $1$-normal and so strict normality is equivalent to $0$-normality. However, the Jacobi Identity coupled with $[T,HM] \subset HM$ implies 
\begin{align*}
\aip{[T,X]}{Y}{} &=  - \aip{ [[T,X],JY]}{T}{} = \aip{ [[X,JY],T]}{T}{} + \aip{ [[JY,T],X]}{T}{} \\
&= T \aip{X}{Y}{} + \aip{ [JY,T]}{JX}{}
\end{align*}
This implies that $0$-normality is equivalent to $\aip{[Y,T]}{X}{} =  -\aip{[T,JY]}{JX}{}$. But this equivalent to  $[T,JY]= J[T,Y]$ which is Tanaka's definition of normal for a strictly pseudoconvex pseudohermitian manifold, \cite{Tanaka}.  
 \enX

The tensors $C^{(j)}$ provide the essential ingredient for the definition of our connections. The idea boils down to using the Levi-Civita connection on each component $\V[j]$ and using projections of  the Lie derivative for mixed components. In general, this will not produce a metric compatible connection, but we can use the tensors $C^{(j)}$ to adjust appropriately.

\bgL{Connection}
If $g$ is an extension of an $r$-graded sRC-manifold, then there exists a unique connection $\nb$ such that
\begin{itemize}
\item $g$ is metric compatible
\item $\V[j]$ is parallel for all $j$
\item $\torr(\V[j],\V[j])  \subseteq \hV[j]$ for all $j$
\item $\aip{ \torr(X^{(j)},Y^{(k)})}{Z^{(j)}}{} = \aip{ \torr(Z^{(j)},Y^{(k)})}{X^{(j)}}{}$ for all $j,k$
\end{itemize}
Furthermore, if $X,Y$ are horizontal vector fields, then $\nb X$ and $\torr(X,Y)$ are independent of the choice of grading and extension $g$. (They do however depend on choice of $VM$)
\enL

\pf For a vector field $Z$, we denote the orthogonal projections of $Z$ to $\V[j]$  by $Z_j$ . Define a new connection $\nb$ as follows: for $X,Y,Z$ sections of $\V[j]$ and $T$ a section of $\hV[j]$ set
\begin{equation}\label{E:CC:Nabla}
\begin{split}
\aip{ \nabla_{X} Y}{ \hV[j]}{}&=0 \\ 
 \aip{\nabla_X Y}{Z}{}  &= \frac{1}{2} \Big( X \aip{Y}{Z}{} + Y\aip{Z}{X}{} -Z\aip{X}{Y}{}  \\ & \quad  -\aip{X}{[Y,Z]}{}- \aip{Y}{[X,Z]}{} +\aip{Z}{[X,Y]}{} \Big)\\    \nabla_T Y & = [T,Y]_j + \frac{1}{2} C^{(j)} (Y,T) \\
  \end{split}
\end{equation}
for $X,Y,Z$ horizontal vector fields and $T,U,W$  vertical vector fields.  It's easy to check that this defines a connection with the desired properties. Futhermore if  $\overline{\nabla}$ is the Levi-Civita connection for $g$, then for sections $X,Y$ of \V[j],
\[ \nabla_X Y = \left( \overline{\nabla}_X Y \right)_j\]
 For uniqueness, suppose that connections $\nabla$ and $\nabla^\prime$ satisfy the required properties and set $A(W,Z) = \nabla_W Z  -\nabla^\prime_W  Z$. Then for sections $X,Y,Z$ of $\V[j]$,  since the torsion terms are in \hV[j] we see 
\begin{align*}
\aip{A(X,Y)}{Z}{} &= - \aip{Y}{A(X,Z)}{} = - \aip{Y}{A(Z,X)}{}\\
&= \aip{A(Z,Y)}{X}{} = \aip{A(Y,Z)}{X}{} \\
&= - \aip{Z}{A(X,Y)}{}
\end{align*}
Similarly if $T$ is a section of $\hV[j]$,
\begin{align*}
\aip{ A(T,X) }{Y}{}&= - \aip{ X}{A(T,Y) }{} = - \aip{ X}{ \tor(T,Y) - \text{Tor}^\prime(T,Y)}{} \\
&= -\aip{\tor(T,X) - \text{Tor}^\prime(T,X)}{Y}{}\\
&= - \aip{A(T,X)}{Y}{}
\end{align*}
Thus $A =0$. Thus this connection $\nabla$ is the unique connection with the desired properties. The required independence from $g$ follows easily from \rfE[CC]{Nabla} and \rfR{indep}
\epf

\bgR{Conn}
An $r$-grading induces a family of connections $\nb[1], \dots ,\nb[r]$ associated to each of sub-gradings of \rfR{subgrading}. Each of these connections agrees in the sense that
\[  \nb[j] X^{(k)}  = \nb X^{(k)} \]
whenever $0 \leq k<j$. In particular, for a horizontal vector fields this means
\[ \nb[1] X = \nb[2] X = \dots = \nb X.\] Thus the differences between the connections can be viewed as a choice in how to differentiate vertical vectors.  The connection $\nb[1]$ associated to the basic grading with be referred to as the basic connection.  We shall denote the basic connection by $\nabla$.
\enR

The following lemma is trivial and left to the reader.

\bgC{Properties}
If $M$ admits an $r$-grading, then 
\begin{itemize}
\item $\torr(\V[j],\V[j])=0$  if and only if $\V[j]$ is integrable.
\item $\torr(HM,\V[j]) \subset HM \oplus \V[j] \oplus \V[j+1]$ for all $j$
\end{itemize}
If the $r$-grading is $j$-normal then
\[ \torr(TM,\V[j]) \subset \hV[j] \]
If the $r$-grading is $0$-normal and  $j$-normal then
\[ 
 \torr(HM,\V[j]) \subseteq \V[j+1] \]
with equality holding if and only if the grading is $j$-regular.
\enC


\bgX{BG}
Suppose that $HM$ has global orthonormal frame $\{X_i\}$ and $VM$ has global orthonormal frame $\{T_\ub\}$ with the following bracket identities:
\begin{align*}
[X_i,X_j] &=c_{ij}^k X_k + c_{ij}^\ua T_\ua  \\
[X_i,T_\ub] &= c_{i\ub}^k X_k + c^\ua_{i \ub} T_\ua \\
[T_\uc,T_\ub] &= c^k_{\uc\ub}X_k + c^\ua_{\uc \ub} T_\ua
\end{align*}
Then using the basic grading and connection we have
\begin{itemize}
\item $VM$ is normal if and only if $c_{i\ub}^k = - c_{k \ub}^i$
\item $g$ is strictly normal  if and only if $c_{i\ub}^\ua = - c_{i \ua}^\ub$ and $c_{i\ub}^k = - c_{k \ub}^i$
\item $g$ is vertically rigid if and only if $\sum c_{i\ub}^\ub =0$ 
\end{itemize}
and 
\begin{itemize}
\item $\nabla_{X_i}X_j = \frac{1}{2} \left( c_{ij}^k +c_{ki}^j + c_{kj}^i \right) X_k $, \qquad $\tor(X_i,X_j) = -c_{ij}^\ua T_\ua$
\item $\nabla_{T_\ub} X_j =  \frac{1}{2} \left( c^j_{k \ub} - c_{j \ub}^k  \right)X_k  $
\item $\nabla_{X_j} T_\ub = \frac{1}{2} \left(c^\ua_{j \ub}  - c^\ub_{j \ua} \right)T_\ua $, 
 \item $\nabla_{T_\uc} T_\ub = \frac{1}{2} \left( c_{\uc \ub}^\ua +c_{\ua \uc}^\ub + c_{\ua \ub}^\uc \right) T_\uc$, \qquad  $\tor(T_\uc,T_\ub) = - c^k_{\uc \ub} X_k$
\item $\tor(X_j,T_\ub) =  -\frac{1}{2} \left( c^j_{k \ub} + c_{j \ub}^k  \right)X_k-\frac{1}{2} \left(c^\ua_{j \ub}  + c^\ub_{j \ua} \right)T_\ua $. 
\end{itemize}
\enX

To illustrate some important behavior, we shall highlight a group particular cases of the previous example

\bgX{C3b}
Let $M$ be the 4 dimensional Carnot group of  \rfX{C3}. Using the basic grading, we can easily compute that 
\[ \nabla_X T = S - \frac{1}{2} S = \frac{1}{2} S, \quad \nabla_X S = 0 - \frac{1}{2} T  = -\frac{1}{2} T,  \]
\[ \tor (X,Y)=-T, \quad \tor(X,T) =-\frac{1}{2} S, \quad \tor(X,S) = - \frac{1}{2} T \]
All other covariant derviatives of frame elements vanish. That the basic covariant derivatives of the natural vertical frame do not vanish is typical of non-step 2 Carnot groups.

However if we use the more refined $2$-grading, then all covariant derivatives of the frame elements vanish and the only non-trivial behavior occurs in the torsion
\[ \torr[2] (X,Y)=-T, \quad \torr[2](X,T) =-S, \quad \torr[2](X,S) =0 \]
\enX

\bgX{sn}
Let $M = \rn{4}$ with the following global orthonormal frames for $HM$ and $VM$
\begin{align*}
X &= \pd{}{x}, \quad Y = \pd{}{y}  + \sin x \pd{}{t} - \cos x \pd{}{s} \\
T &= \cos x \pd{}{t} + \sin x \pd{}{s}   , \quad S=-\sin x \pd{}{t} + \cos x \pd{}{s}
\end{align*}
Then $[X,Y]=T= -[X,S]$, $[X,T]=S$ with all other commutators vanishing. It's then easy to check that this is  a strictly normal extension for the basic grading and that the only non-trivial covariant derivatives are then $\nabla_X T = S$ and $\nabla_X S =T$.  This is an example of a flat, equiregular,  strictly normal sRC-manifold with step size $>2$.
\enX

\bgX{spsic2}
Let $(M,J,\eta)$ be a strictly pseudoconvex pseudohermitian manifold. The Tanaka-Webster connection is the unique connection such that $\eta, d\eta$ and $J$ are parallel and the torsion satisfies
\[ \tor(X,Y) = d\eta(X,Y)T, \qquad \tor(T,JX) = - J \tor(T,X) \]
The only defining property of the basic connection not clearly satisfied by the Tanaka-Webster connection is torsion symmetry. But if we pick $X,Y$ as any horizontal vector fields then the Jacobi identity implies
\begin{align*}
0 &=\eta\left(  [ T,[X,JY]]+[JY,[T,X]] + [X,[JY,T]] \right)\\
&= -T\aip{X}{Y}{} + \aip{[T,X]}{Y}{} +\aip{[JY,T]}{JX}{}\\
&= -\aip{X}{\nabla_T Y}{} +\aip{\tor(T,X)}{Y}{} -\aip{\nabla_T JY}{JX}{} + \aip{\tor(T,JY)}{JX}{}\\
&= \aip{\tor(T,X)}{Y}{} + \aip{\tor(T,JY)}{JX}{}\\
&= \aip{\tor(T,X)}{Y}{}-\aip{\tor(T,Y)}{X}{}.
\end{align*}
Thus the Tanaka-Webster connection satisfies the requirements of the basic connection.
\enX

One of the key computational tools when using the Levi-Civita connection is the existence of Riemannian normal coordinates in a neighborhood of any given point. As $HM$ is non-integrable in every interesting example, we cannot expect to find a similarly useful coordinate system in the subRiemannian case. However, when the extension is normal, we can guarantee the existence of a local orthonormal horizontal frame with computationally nice properties at any particular point $p$.

\bgD{NNF}
If $M$ is $r$-graded, then an orthonormal frame $\{E^{(j)}_i\}$ for $\V[j]$, $0 \leq j \leq r$, defined in a neighborhood of $p$ is  $\nb$-normal at $p$ if
\[ \left( \nb  E^{(j)}_i \right)_{| p} = 0\]
\enD

\bgL{NormalC}
Suppose $g$ is a $j$-normal $r$-grading. Then there exists a $\nb$-normal frame for $\V[j]$ at every $p\in M$.
\enL

\pf
Let $v^{(k)}_1,\dots,v^{(k)}_{n_k}$ be orthonormal vectors spanning $\V[k]_p$ and let $\{x_{(k)}^i\}$ be the coordinates near $p$ induced by the exponential map of $\nb$ at $p$ using this frame.  Then certainly 
$ \nb_{c_i \partial_{x^i}} ( c_i\partial_{x^i} ) =0$ at $p$ whenever the coefficients $c_i$ are constant. Considering $\nb_{ \partial_i + \partial_j} ( \partial_i + \partial_j) $ in particular, this implies that for all $i,j$ at $p$
\[ 0 = \nb_{\partial_i} \partial_j + \nb_{\partial_j} \partial_i  = 2 \nabla_{\partial_i} \partial_j + \torr(\partial_i,\partial_j) .\]
Now $(\partial_{x^i_{(j)}})_{|p} \in \V[j]_p$. Since torsion is tensorial and $\torr(TM,\V[j])\subset \hV[j]$ by \rfC{Properties}, this implies that 
\bgE{nj}
\left( \nb_{\partial_i} \partial_{x^i_{(j)}} \right)_{|p} \in \hV[j]_p
\enE
for all $i$.

Now in a small neighbourhood of $p$ define $Z^{(j)}_k = ( \partial_{x^{(j)}}^k)_0$, i.e. the orthogonal projection of $ \partial_{x^{(j)}}^k$ onto $\V[j]$. Set $T^{(j)}_k = \partial_{x_{(j)}} - Z^{(j)}_k$.  We clearly have linear independence near $p$ and so  $Z^{(j)}_1,\dots,Z^{(j)}_{n_j}$ is a local frame for $\V[j]$.

Now for any vector field $Y$,
\[ \nb_{Y}Z_i^{(j)} = \nb_{Y} ( \partial_{x^{(j)}}^i - T_i^{(j)}) =   \nb_{Y}  \partial_{x^{(j)}}^i - \nb_{Y} T_i^{(j)}\]
The first  term on the right is in \hV[j] by \rfE{nj}. The last term is in \hV[j] everywhere as $T_i^{(j)}$ is a section of $\hV[j]$ which is parallel. But  $\nabla_Y Z^{(j)}_i$ is in \V[j] as $\V[j]$ is parallel. This implies that $\nb Z_i^{(j)}=0$ at $p$.

Now from metric compatibility, we see that  $Y \aip{Z^{(j)}_i}{Z^{(j)}_k}{|p} = 0$ for each $i,k$, so an easy induction argument shows that if we apply the Gram-Schmidt algorithm to $Z^{(j)}_1,\dots, Z^{(j)}_{n_j}$ we obtain an orthonormal frame with the same property at $p$.

\epf

\bgC{NormalC}
If the grading is strictly normal,  then near any point $p \in M$, there is a graded orthonormal frame $X^{(j)}_i$ for $TM$ such that $\left( \nabla  X^{(j)}_i \right)_{| p} = 0$.
\enC

\section{Curvature and the Bianchi Identities}\setS{CV}

\bgD{CV}
The subRiemannian curvature  tensors for a sRC-manifold with extension $g$ are defined by
\[ R(A,B)C = \nabla_A \nabla_B  C - \nabla_B \nabla_A C - \nabla_{[A,B]} C \]
and
\[  \srm(A,B,C,D) = \aip{ R(A,B)C }{D}{} \]
\enD

We note that for any vectors $A,B \in TM$, the restriction of the $(1,1)$-tensor $R(A,B)$ to $HM$ is independent of the choice of extension $g$.

This definition immediately yields notions of flatness in subRiemannian geometry.

\bgD{Flat}
We say that an $M$ is horizontal flat if  $\srm(\cdot,\cdot, HM, \cdot)=0$ for any extension $g$. A particular extension is vertically flat if $\srm(\cdot,\cdot, VM, \cdot)=0$ or flat if $\srm =0$.
\enD

\bgL{HFlat} 
A sRC-manifold is horizontally flat if and only if in a neighborhood of every point $p\in M$ there is a local orthonormal frame $\{E_i\}$ for $HM$ such that $\nabla E_i =0$. If $HM$ is integrable, this local frame can be chosen to be a coordinate frame.

A similar result holds for a vertically flat extension $g$ and $VM$.
\enL

The proof of this lemma follows that of Theorem 7.3 in \cite{Lee3} almost exactly. 

\bgX{CarnotF}
Every step $r$ Carnot group is horizontally flat for the basic grading and flat for the $r-1$-grading.  The sRC-manifolds considered in \rfX[CC]{C3} and \rfX[CC]{sn} are both flat.
\enX

For convenience of notation, it is useful to define the following

\bgD{Cyclic}
If $S$ is any set and $F\colon S^k \to L$ is any map into a vector space $L$, we define  $\C F$ to be the sum of all cyclic permutations of $F$. For example if $k=3$, then
\[ \C F (X,Y,Z) = F(X,Y,Z) + F(Y,Z,X) + F(Z,X,Y)\]
\enD

An example of the cyclic construction in action is a compressed form of the Jacobi Identity for vector fields, namely 
\[ \C [X,[Y,Z]]   =0 \]
We shall use it primarily to efficiently describe symmetries of the curvature tensor.





We also introduce 
\bgD{CT}
The second-order torsion of $\nabla$ is the $(3,1)$-tensor 
\[ \Tor(A,B,C) = \tor(A,\tor(B,C)) \]
\enD

We are now in a position to discuss the fundamental questions of curvature symmetries. Many of the properties of the Riemannian curvature tensor go through unchanged, with exactly the same proof. In particular,
\bgL{Symmetries}
The subRiemannian curvature tensor always has the following symmetries
\begin{itemize}
\item $\srm(A,B,C,D) = -\srm(A,B,D,C)$
\item $\srm(A,B,C,D) =-\srm(B,A,C,D)$
\item $\srm(TM,TM,HM,VM) = 0 $
\end{itemize}
\enL
However, many symmetry properties of the Riemannian curvature tensor require additional assumptions in the subRiemannian case. Most of these symmetries are naturally related to the Bianchi Identities. 

\bgL[Algebraic Bianchi Identites]{HAB}
For any vector fields $X,Y,Z,$
\[ \C R(X,Y)Z = -\srt(X,Y,Z) + \C (\nabla \tor)(X,Y,Z).\]
Furthermore 
\begin{enumerate}
\item if $X,Y,Z \in \V[j]$ then
  \[ \C (\nabla \tor)(X,Y,Z) \in \hV[j] \]
 \item if $X,Y,Z \in \V[j]$ and the grading is $j$-normal, then
  \[ -\srt (X,Y,Z) \in \hV[j] \]
 \item if $X,Y \in \V[j]$, the grading is $j$-normal and $\hV[j]$ is integrable then
 \[ -\srt (X,Y,Z) \in \hV[j] \]
  \end{enumerate}
\enL

\pf The first part of the lemma is a standard result from differential geometry, but for completeness we shall present a short proof
\begin{align*}
 \C R (X,Y)Z &= \C \left( \nabla_X \nabla_Y  Z - \nabla_Y \nabla_X Z - \nabla_{[X,Y]}Z \right) \\
 &= \C  \left( \nabla_Z \nabla_X Y - \nabla_Z \nabla_Y X - \nabla_{[X,Y]}Z \right)\\
 &= \C \left( \nabla_Z ( [X,Y]+ \tor(X,Y) ) - \nabla_{[X,Y]} Z \right) \\
 &= \C \left( [Z,[X,Y]] + \tor(Z,[X,Y]) + \nabla_Z \tor(X,Y) \right)\\
 &= \C \left( \tor(Z,[X,Y])  +\tor(\nabla_Z X,Y)  + \tor(X,\nabla_Z Y) \right)  \\ & \qquad+ \C (\nabla \tor)(X,Y,Z) \\
&= \C \left(\tor(Z, [X,Y] - \nabla_X Y - \nabla_Y X\right)) + \C (\nabla \tor)(X,Y,Z)\\
&= -\srt(X,Y,Z) + \C (\nabla \tor)(X,Y,Z)
\end{align*}

The remaining parts consist of analyzing the terms $\srt$ and $\C (\nabla \tor)$. Since these are tensorial, we can compute using normal and seminormal frames. First let $X,Y,Z$ be elements of a seminormal frame for $\V[j]$ at $p$, then
\[
 \C (\nabla \tor) (X,Y,Z) = \C \left( \nabla_X \tor(Y,Z) \right) 
\]
But each torsion piece must be in $\hV[j]$. As this is bundle parallel,  we have established (a).  Now, if we assume the frame is $j$-normal, then we can instead use a normal frame at $p$. If $X$ is an element of this frame then $\tor(X,TM) \subset  \hV[j]$ and it is easy to check that (b) holds.

Finally, assume a $j$-normal grading and that $X,Y$ are elements of a $j$-normal frame at $p$, but $Z$ is an arbitrary vector field.  then
\[ -\srt(X,Y,Z) = \C \tor(X,\tor(Y,Z)) = \tor(Z,\tor(X,Y) ) \]
But 
\begin{align*}
\left( \tor(Z,\tor(X,Y) ) \right)_j = - [Z_{\hat{j}},\tor(X,Y)]_j  
\end{align*}
which vanishes if $\hV[j]$ is integrable. Thus (c) holds.

\epf

\bgC[Horizontal Algebraic Bianchi Identity]{HAB}
If $X,Y,Z,W$ are horizontal vector fields and $VM$ is normal, then
\[ \aip{ \C R(X,Y)Z }{W} {} = 0\]
If $VM$ is also  integrable, then this can be relaxed to any three of $X,Y,Z,W$ horizontal. 
\enC

\bgC{AHBc}
If $VM$ is normal and  $X,Y,Z,W$ are horizontal vector fields then
\[ \srm (X,Y,Z,W) =\srm (Z,W,X,Y )\]
If $VM$ is also  integrable, then this can be relaxed to any three of $X,Y,Z,W$ horizontal. 
\enC

The proof is nearly identical to the Riemannian case, see \cite{Lee3}, Proposition 7.4 for example. We shall only briefly sketch out the details.

\pf A straightforward computation shows that
\bgE{AHBc}  2 \srm(C,A,B,D)-2 \srm(B,D,C,A) =  \C \aip{ \C R(A,B)C }{D}{} \enE
The result then follows from the horizontal algebraic Bianchi Identity,

\epf

\bgL[Differential Horizontal  Bianchi Identites]{DHB}
For any vector fields $X,Y,Z,V$
\[ \C \left(\nabla_W R(X,Y) \right) Z = \C \left( R(\tor(X,V), Y) \right)Z \]
Furthermore, if  $VM$ is normal and integrable and $X,Y,Z,W,V \in HM$ then
\[ \nabla \srm(X,Y,Z,W,V) + \nabla \srm(X,Y,W,V,Z) + \nabla \srm(X,Y,V,Z,W) = 0 \]
\enL

\pf  Again, the first part is a standard result that can be derived as follows
\begin{align*}
\left( \nabla_WR\right) (X,Y)Z &= \nabla_WR(X,Y)Z - R(\nabla_WX,Y)Z - R(X,\nabla_WY)Z \\ & \qquad - R(X,Y)\nabla_WZ  \\ 
&= [ \nabla_V, R(X,Y)] Z - R(\nabla_WX,Y)Z - R(X,\nabla_WY)Z 
\end{align*}
Thus, recalling the Jacobi identity applies to operators, we see
\begin{align*}
 \C \left( (\nabla_W) R(X,Y) \right) Z &= \C \left( [ \nabla_W,R(X,Y)] \right) Z  \\& \qquad - \C \left( R(\nabla_W X,Y) \right) Z  - \C\left(R(X,\nabla_WY) \right)Z \\
 &= \C \left( [ \nabla_V, [\nabla_X,\nabla_Y] ] - [\nabla_V,\nabla_{[X,Y]} ] \right) Z \\
 & \qquad - \C \left( R(\nabla_WX,Y) \right) Z  +  \C\left(R(\nabla_X V,Y) \right)Z \\
 &=  - \C \left( [\nabla_V,\nabla_{[X,Y]} ] \right) Z+  \C \left( R([X,V]+\tor(X,V), Y) \right)Z\\
 &=  - \C \left( [\nabla_V,\nabla_{[X,Y]} ] \right) Z+  \C \left( R(\tor(X,V), Y) \right)Z\\
 &   \qquad +\C \left(  [\nabla_{[X,V]},\nabla_Y] -\nabla_{[[X,V],Y]} \right)Z\\
  &=  - \C \left( [\nabla_V,\nabla_{[X,Y]} ] \right) Z+  \C \left( R(\tor(X,V), Y) \right)Z\\
  & \qquad +\C \left( [ \nabla_{[Y,X]}, \nabla_V] \right)Z \\
  &=  \C \left( R(\tor(X,V), Y) \right)Z
 \end{align*}
 To see the second part, we note that as $VM$ is normal \rfC{AHBc} implies that the required identity is equivalent to 
 \[ \nabla \srm(Z,W,X,Y,V) + \nabla \srm(W,V,X,Y,W,Z) + \nabla \srm(V,Z,X,YW) = 0 \]
Choose $X,Y,Z,W,V$ to be elements of a normal from for $HM =\V[0]$ at $p$, then
\begin{align*}
\nabla \srm(Z,W,X,Y,V) & + \nabla \srm(W,V,X,Y,W,Z) + \nabla \srm(V,Z,X,YW) \\
&= \aip{\C \left( (\nabla_V) R(Z,W) \right) X}{Y}{} \\
&= \aip{ \C \left( R(\tor(Z,V), W) \right)X}{Y}{}
\end{align*}
 But by \rfC{AHBc}
 \[  \srm( \tor(Z,V),W,X,Y) = \srm(X,Y,\tor(Z,V),W) = 0 \]
 as $\tor(Z,W)$ is vertical.
\epf

\section{Ricci Curvature and Bochner Formulae}\setS{RC}

\bgD{Ricci}
We define the subRiemannian Ricci curvature of $\nabla$ by
\[ \Ric(A,B) = \sum\limits_k  \srm(A,X_k,X_k,B) \]
where $\{X_k\}$ is any horizontal orthonormal frame. The horizontal scalar curvature is defined by
\[ S_0 = \tr \Ric =  \Ric(X_k,X_k)\]
It should be noted that the scalar curvature is independent of the choice of extension $g$ as is the Ricci curvature restricted to horizontal vector fields.
\enD



It should be remarked here, that in general the Ricci curvature for the canonical connection is not symmetric. However, using \rfC[CV]{AHBc} and elementary properties of the connection, we can immediately deduce 
\bgL{Ricci}
If $VM$ is normal and $X,Y \in HM$ then 
\[ \Ric(X,Y) = \Ric(Y,X) \]
If  $VM$ is normal and integrable then
\[ \Ric(VM,HM)=0\]
\enL

\pf
The first  follows from the corollary to the horizontal Bianchi Identity.  For the second, we apply the corollary to the horizontal Bianchi Identity, to see that
\[ \Ric(U,X) = \srm(E_k,U,X,E_k) = \srm(X,E_k,E_k,U) =0.\]
\epf

\bgL[Contracted Bianchi Identity]{CBI} Suppose $VM$ is normal and integrable, then for any horizontal $X$
\[ \nabla_X S_0 = 2 \sum (\nabla  \Ric)(E_j,X,E_j)\]
where $E_i$ is an orthonormal frame for $HM$. Equivalently
\[ \nabla_0 S_0 = 2 \tr (\nabla \Ric) \]
\enL

\pf Let $X$ be any element of a normal frame at $p$. Apply the differential Bianchi Identity to $E_i,E_j,E_j,E_i,X$ and sum over $i$ and $j$.

\epf

As a quick and easy consequence of this identity, we get a subRiemannian version of a result of Schur, that whenever the Ricci tensor is conformally equivalent to the metric then the manifold is Einstein.

\bgC{Schur} Suppose that $M$ is a connected sRC-manifold such that  $HM$ bracket generates, $\text{dim}(HM)=d >2$ and that $VM$ is normal and integrable. If
\[ \Ric(X,Y) = \lambda \aip{X}{Y}{}\]
for horizontal all vectors $X,Y$ then $\lambda$ must be constant.
\enC

\pf Let $E_i$ be a normal frame at $p \in M$. Then at $p$,
\[ S_0 = \Ric(E_i,E_i) = \lambda d \]
but
\[ 2 \tr (\nabla \Ric)(E_j) = 2 \nabla_i \Ric (E_j,E_i) = 2 E_j \lambda \]
Since $E_j S_0 = 2 \tr (\nabla \Ric) (E_j)$, we must have $d=2$ or $E_j \lambda =0$. Thus all horizontal vector fields annihilate $\lambda$. As $HM$ bracket generates, this implies that $\lambda$ is constant 

\epf

One of our purposes is to use Bochner type results to study the relationships between curvature, geometry and topology on subRiemannian manifolds. To use this theory, we shall need a geometrically defined subelliptic Laplacian.

\bgD{HL}
For a tensor $\tau$, the horizontal gradient of $\tau$ is defined by
\[ \nabla_0 \tau = \nabla_{E_i} \tau \otimes E_i, \]
where $E_i^*$ is the dual to $E_i$.

The horizontal Hessian of $\tau$ is defined  by
\[ \nabla^2_0 \tau (X,Y)=\left(  \nabla_{X} \nabla_{Y} - \nabla_{\nabla_{X}Y} \right)  \tau\]
for $X,Y \in HM$ and zero otherwise.

The symmetric horizontal Hessian of $\tau$ is defined by
\[ \nabla^{2,sym}_0 \tau(X,Y) = \frac{1}{2} \left(\nabla^2_0 \tau (X,Y) +\nabla^2_0 \tau (Y,X) \right)\]
Finally, the horizontal Laplacian of $\tau$ is defined by
\[ \lp \tau = \tr \left( \nabla^2_0 \tau \right)  =\left(  \nabla_{E_i} \nabla_{E_i} - \nabla_{\nabla_{E_i}E_i} \right)  \tau  \]
\enD

The Laplacian on a Riemannian manifold has a rich and interesting $L^2$-theory. To replicate this for sRC-manifolds, it is necessary to choose a metric extension. This metric extension then yields a volume form and we have meaningful $L^2$-adjoints. Unfortunately, the horizontal Laplacian defined here, does not always behave as nicely as the Riemannian operator. However, if we make a mild assumption on the metric extension, much of the theory can be generalized. 

\bgD{VR}
For a metric extension of an $r$-grading we define a $1$-form $\mathfrak{R}_g$ by
\[\mathfrak{R}_g (v)= \sum\limits_{j>0} \sum\limits_i   B^{(j)} (E_i^{(j)},E_i^{(j)},v_0) \] where $E^{(j)}_i$ is an orthonormal frame for $\V[j]$ .

We say that a complement $VM$ is vertically rigid if there exists a metric extension $g$ such that
\[ \mathfrak{R}_g \equiv 0.\]
\enD

\bgL{VReq}
For an orientable sRC-manifold, the following are equivalent
\begin{enumerate}
\item $VM$ is vertically rigid
\item There exists a volume form $dV$ on $M$ such that for any horizontal vector field
\[ \text{div } X = \text{tr}_0 \nabla X = \aip{\nabla_{e_i} X}{e_i}{}\]
where $e_i$ is an orthonormal frame for $HM$. 
\item Every metric extension $g$ is vertically conformal to a metric $\tilde{g}$ with $\mathfrak{R}_{\tilde{g}} \equiv 0$
\end{enumerate}
Furthermore, if $HM$ bracket generates, then the volume form in (b) is unique up to constant multiplication.
\enL

\pf 

To show that (a) implies (b), we first note that for the particular metric extension  $g$ with $\mathfrak{R}_g \equiv 0$, we have
 \[ \sum\limits_{j>0} \sum\limits_i   \aip{ \tor(E^{(j)}_i,X)}{E^{(j)}_i}{}  = 0.\]
Now we recall the standard result  (see \cite{Kobayashi}, appendix 6) that since $\nabla$ is metric compatible and $HM$ is parallel, the divergence operator for the metric volume form $g$ satisfies
 \[\begin{split}  \text{div}_g X  &= \text{tr} (\nabla +\tor )(X) \\
&=\text{tr}_0 \nabla X + \sum\limits_{j>0} \sum\limits_i   \aip{ \tor(E^{(j)}_i,X)}{E^{(j)}_i}{}  \\
&=\text{tr}_0 \nabla X - \mathfrak{R}_g (X)\\
&=\text{tr}_0 \nabla X
\end{split}\]
Thus we can set $dV=dV_g$.

To show (b) implies (c), we consider metrics vertically conformal to an arbitrary extension $g$,
\[ g_{\lambda}  = \begin{cases}  g, \quad \text{on $HM$}\\ e^{\lambda} g, \quad \text{on $VM$} \end{cases}.\]
Now if $dV_g = e^{\mu} dV$, then set $\lambda = -\frac{\mu}{\text{dim}(VM) }$ so $dV_{g_\lambda} = dV$.  Then for horizontal $X$
\[ \text{tr}_0 \nabla X - \mathfrak{R}_{g_{\lambda}} (X) = \text{div}_{g_\lambda} X = \text{div }X =\text{tr}_0 \nabla X  \]
so $\mathfrak{R}_{g_{\lambda}} \equiv 0$. 

Since (c) trivially implies (a), the equivalence portion of the proof is complete. For the uniqueness portion, we note that if $\Omega =e^{\lambda }dV$ then for any horizontal $X$, we have
\[ \text{div}_{\Omega} X = \text{div } X - X(\lambda) .\]
If the two divergences agree on horizontal vector fields and $HM$ bracket generates, this immediately implies that $\lambda$ is a constant.

\epf

For an orientable, vertically rigid sRC-manifold, there is then a $1$-dimensional family of volume forms for which $\text{div }X =\text{tr}_0 \nabla X$. We shall often refer to such a volume form as a rigid volume form. Vertical rigidity therefore gives us a canonical notion of integration on a sRC-manifold that does not depend on the choice of metric extension.
 
 As an immediate consequence, we have

\bgL{VR}
Suppose that $M$ is orientable and $VM$ is vertically rigid.  Then on functions,
\[ \lp  = E_i^2 + \text{div } E_i =- \nabla^*_0 \nabla_0  \]
where the divergence and $L^2$ adjoint are taken with respect to a rigid volume form.
\enL
Thus on a vertically rigid sRC-manifold, the horizontal Laplacian behaves qualitatively in a similar fashion to the Riemannian Laplacian.

\bgR{VR}
A slightly different definition of vertically rigid was presented in \cite{HP2}, \cite{HP4}. All the results of those papers that required vertical rigidity, can be obtained using this definition with only very minor modifications to the proofs.
\enR

Our main result utilizing the horizontal Laplacian is the following subRiemannian equivalent of the classical Bochner formulas.

\bgT{Bochner} If $F$ is a closed vector field and $F_k$ is the projection of $F$ to $\V[k]$ then 
\[ \begin{split} \frac{1}{2} \lp  \left| F_j \right|^2  &=\Ric(F_j,F_0)  + \left| \nabla_{(0)} F_j \right|^2\\ 
& \qquad  + \sum\limits_{i} \Big(  \aip{E_i}{\nabla^2 F_j ( F_j, E_i) }{} \\
& \qquad -  2\aip{\nabla_{E_i} F}{\tor(E_i,F_j)}{}  + \aip{F}{ (\nabla \tor) (F_j,E_i,E_i) }{}  \\
& \qquad - \aip{F}{\Tor(E_i,E_i,F_j)}{} \Big) \\
\end{split}\]
where $\{E_i\}$ is any orthonormal horizontal frame. 
\enT

Before we prove this result, we introduce some terms and notation. Firstly, we define $J \colon TM \times TM \to TM$ by
 \bgE{J} \aip{J(A,Z)}{B}{} =\aip{ \tor(A,B)}{Z}{}\enE
Next we recall that a vector field $F$ is closed if 
 \[ A \mapsto \aip{F}{A}{} \]
 is a  closed $1$-form. It is then easy to check that $F$ is closed if and only if for all vector fields $A,B$
\begin{align*}
 \aip{\nabla_B F}{A}{}&=\aip{ \nabla_A F}{B}{} - \aip{J( B,F)}{A}{} = \aip{ \nabla_A F}{B}{} + \aip{J(A,F)}{B}{} 
\end{align*}
 
\begin{proof}
Set  $u = \frac{1}{2} \left| F_j \right|^2$,  then
\bgE{use1}
\begin{split}
\aip{\nabla_{(0)} u}{Y}{} &= \aip{ \nabla_Y F_j}{F_j}{}=\aip{ \nabla_Y F}{F_j}{}  \\
&= \aip{\nabla_{F_j} F_0}{Y}{} + \aip{J(F_j,F) }{Y}{}  
\end{split}
\enE
so
\[ \nabla_{(0)} u =  \nabla_{F_j} F_0 + J( F_j,F)_0.\]
Next we need some preliminaries. Firstly, for horizontal $X,Y$ 
\bgE{use2}
\begin{split}
\nabla^2 u (X,Y)& = X \aip{Y}{\nabla u}{} - \aip{\nabla_X Y}{\nabla u}{} = \aip{Y}{\nabla_X \nabla_0 u}{} \\\ & = \aip{Y}{\nabla_X \nabla_{(0)} u}{}
\end{split}
\enE
Secondly,
\begin{align*}
\aip{\nabla_X J(F_j,F)}{X}{} &= X \aip{ J(F_j,F)}{X}{} - \aip{J(F_j,F)}{\nabla_X X}{}\\
&= X \aip{F}{\tor(F_j,X)}{} - \aip{J(F_j,F)}{\nabla_X X}{} \\
&= \aip{\nabla_X F}{\tor(F_j,X) }{} + \aip{F}{\nabla_X \tor(F_j,X)}{} \\& \qquad  - \aip{F}{\tor(F_j,\nabla_X X)}{} \\
\end{align*}

Now we can begin the main computation. For horizontal $X$
\begin{align*}
\nabla^2 u(X,X) &= \aip{\nabla_X \nabla_{(0)} u}{X}{} = \aip{ \nabla_X \nabla_{F_j} F_0}{X}{} + \aip{\nabla_X J(F_j,F)}{X}{} \\
&= R(X,F_j,F_0,X) + \aip{\nabla_{F_j} \nabla_X F_0}{X}{} + \aip{\nabla_{[X,F_j]}F_0 }{X}{} \\ & \qquad + \aip{\nabla_X J(F_j,F)}{X}{} \\
&= R(X,F_j,F_0,X)  +  \aip{\nabla_{F_j} \nabla_X F_0}{X}{} \\
& \qquad + \aip{\nabla_{\nabla_X F_j -\nabla_{F_j} X - \tor(X,F_j)} F }{X}{} + \aip{\nabla_X J(F_j,F)}{X}{} \\
&= R(X,F_j,F_0,X)  +  \aip{\nabla_{F_j} \nabla_X F_0  - \nabla_{\nabla_{F_j} X} F_0}{X}{} \\
& \qquad + \aip{\nabla_X F_j}{\nabla_X F_j}{} + \aip{F}{\tor(X,\nabla_X F_j)}{} \\
& \qquad -\aip{\nabla_X F}{\tor(X,F_j)}{} - \aip{F}{\tor(X,\tor(X,F_j))}{}  \\ & \qquad + \aip{\nabla_X J(F_j,F)}{X}{} \\
&= R(X,F_j,F_0,X)  + \aip{X}{\nabla^2  F_0 ( F_j, X) }{} + \left| \nabla_X F_j \right|^2\\
& \qquad -  2\aip{\nabla_X F}{\tor(X,F_j)}{}  + \aip{F}{ (\nabla \tor) (F_j,X,X) }{}  \\
& \qquad - \aip{F}{\Tor(X,X,F_j)}{}
\end{align*}
Now we let $X$ range over the frame $E_i$ and take a sum.
\end{proof}

To apply this theorem, we make the following observations

\bgL{int}
If $F =\nabla f$ then
\[\sum\limits_{i} \big(  \aip{E_i}{\nabla^2 F_0 ( F_j, E_i) }{} \big) = \aip{\nabla_{(j)} f}{\nabla_{(j)} \lp f}{}\]
\enL 

\pf

By \rfE{use2},
\begin{align*}
(\nabla_Z \nabla^2 f) (X,X) &= Z \aip{X}{\nabla_X F_0}{} - \aip{\nabla_Z X}{\nabla_X F_0}{} -\aip{X}{\nabla_{\nabla_Z X} F_0}{}\\
&= \aip{X}{\nabla_Z \nabla_X F_0}{}- \aip{X}{\nabla_{\nabla_Z X} F_0}{}\\
&= \aip{X}{ \nabla^2 F_0 (Z,X)}{}
\end{align*}
so
\begin{align*}
 \sum\limits_{i} \big(  &\aip{E_i}{\nabla^2 F_0 ( F_j, E_i) }{} \big) =  \sum\limits_{i} (\nabla_{F_j} \nabla^2 f)(E_i,E_i) \\
&= \sum\limits_{k} \aip{\nabla f}{U^{(j)}_k }{} \Big( \nabla_{U^{(j)}_k} \lp f \\ & \qquad  - \sum\limits_{i,m} \aip{\nabla_{U^{(j)}_k} E_i}{E_m}{}  \left( \nabla^2 f( E_m,E_i) +\nabla^2 f(E_i,E_m) \right)\Big )\\
&= \aip{\nabla_{(j)} f}{\nabla_{(j)} \lp f}{}
 \end{align*}
 as the latter term is skew-symetric in $i$ and $m$.

\epf

\bgD{BGT}
The Baudoin-Garofalo tensor for an sRC-manifold with metric extension is the unique symmetric 2-tensor such that
\bgE{BGT}
\begin{split}  \mathcal{R}(A,A) & = \Ric(A_0,A_0) + \aip{A}{\tr (\nabla \tor)(A_0)}{} \\ & \qquad +  \frac{1}{4} \sum\limits_{i,j} \left| \aip{\tor(E_i,E_j)}{A}{}\right|^2
\end{split}\enE
\enD

Note that from standard polarization arguments this defines
\[ \mathcal{R}(A,B) = \frac{1}{4} \left( \mathcal{R}(A+B,A+B) - \mathcal{R}(A-B,A-B)  \right)\]

\bgC{Bochnerf}
If $g$ is strictly normal  with respect to the basic grading and $VM$ is integrable then
\[ \begin{split} \frac{1}{2} \lp \left| \nabla_{(0)} f\right|^2 &- \aip{\nabla_{(0)} f}{\nabla_{(0)} \lp f}{} =\\ &\mathcal{R}(\nabla f,\nabla f)  + \| \nabla^{2,sym}_{(0)} f \|^2 -  2\aip{\nabla_{E_i} \nabla_{(1)} f }{\tor(E_i, \nabla_{(0)}f )}{} \\
\frac{1}{2} \lp \left| \nabla_{(1)} f\right|^2 &- \aip{\nabla_{(1)} f}{\nabla_{(1)} \lp f}{} =\| \nabla_{(0)} \nabla_{(1)}  f \|^2 \end{split}\]
\enC

\pf Most of this result follows immediately from noticing that the strictly normal condition eliminates many of the torsion terms  from \rfT{Bochner} and then applying \rfL{int}. The rest consists of analyzing the $\|\nabla_{(0)}^2 f\|^2$. First note
\begin{align*}
\nabla^2 f( E_i,E_j)  &  = E_i E_j f  - (\nabla_{E_i}E_j) f   \\
&= \frac{1}{2} \left( E_iE_jf - (\nabla_{E_i}E_j) f  \right)  +\frac{1}{2} \left( E_jE_i f - (\nabla_{E_j}E_i )f \right) \\ & \qquad + \frac{1}{2} \left( [E_i,E_j]f + (\nabla_{E_j}E_i) f - (\nabla_{E_i}E_j) f \right)\\
&= \frac{1}{2} \left( \nabla^2 f( E_i,E_j)+ \nabla^2 f( E_j,E_i) \right) - \frac{1}{2} \tor(E_i,E_j)f 
\end{align*}
From this we immediately obtain,
\[ \| \nabla^2_{(0)} f\|^2 = \| \nabla^{2,sym}_{(0)} f\|^2 + \frac{1}{4}  \sum\limits_{i,j} \left| \aip{\tor(E_i,E_j)}{\nabla f}{} \right|^2\]
\epf
\
\bgD{tb} The torsion bounds of $M$ are the defined by
\[ \kappa_{ij}^m =  \sup \left\{  \big| \tor(X^{(i)},X^{(j)})_m\big|^2 \colon  \big|X^{(i)} \big|,\big|X^{(j)} \big|\leq1 \right\} \]
Noting that $0 \leq \kappa_{ij}^m \leq +\infty$.
\enD
To obtain topological and geometric information from this result, we follow the technique developed by Baudoin and Garofalo in  \cite{BaudoinGarofalo}. We define symmetric bilinear forms by
\begin{align*}
\Gamma_{(j)} (f,g) &= \aip{\nabla_{(j)} f}{\nabla_{(j)}g }{} \\
\Gamma_{(j)}^2 (f,g) &= \lp \Gamma_{(j)} (f,g) - \Gamma_{(j)} (\lp f,g) - \Gamma_{(j)} ( f, \lp g)
\end{align*}

If $g$ is strictly normal then it is easy to check that
\bgE{com} \Gamma_{(0)} \left(f, \Gamma_{(1)} (f,f) \right) =  \Gamma_{(1)} \left(f, \Gamma_{(0)} (f,f) \right)\enE
and we obtain the following result
\bgT{gci}
Suppose $g$ is strictly normal for the basic grading and $VM$ is integrable. If  
\[  \kappa_{00}^1 < \infty \] 
and there exist constants $\rho_1 \in \rn{}$ and  $\rho_2 >0$ such that\begin{align*}
\mathcal{R}(A,A) &\geq \rho_1 \|A_0\|^2 + \rho_2 \|A_1\|^2  \\
\end{align*}
then for $\kappa = \dim(HM) \kappa_{00}^1$,  the generalized curvature-dimension inequality
\[  \Gamma^2_{(0)} + \nu \Gamma^2_{(1)} \geq \frac{1}{\dim HM} \left(\lp f\right)^2 +\left( \rho_1 - \frac{\kappa}{\nu} \right) \Gamma_{(0)} (f,f)  + \rho_2 \Gamma_{(1)}(f,f)\]
holds for every $f \in C^\infty(M)$ and $\nu>0$ .
\enT

\pf As 
\[ \| \nabla^{2,sym}_{(0)} f\|^2  \geq  \sum\limits_i \left(  \nabla^2f (E_i,E_i)   \right)^2 \geq \frac{1}{\dim HM}  \left( \lp f\right)^2\]
this follows immediately from \rfC{Bochnerf} and the elementary identity that
\[ 2 \sum\limits_i  \left| \aip{\nabla_{E_i} \nabla_{(1)} f }{\tor(E_i, \nabla_{(0)}f )}{} \right|  \leq \nu \| \nabla_{(0)}\nabla_{(1)} f\|^2  + \frac{\kappa}{\nu} \|\nabla_{(0)}f \|^2\]
\epf

It was shown in \cite{BaudoinGarofalo}, that under the additional mild hypothesis that there exists a sequence $h_k \in C_c^\infty(M)$ of increasing functions that converge pointwise to $1$ everywhere and satisfy
\[  \left\| \Gamma_{(0)}(h_k,h_k) \right\|_\infty + \left\|\Gamma_{(1)}(h_k,h_k)\right\|_\infty \to 0,\]
the generalized curvature inequality has a wide variety of topological, geometric and analytical consequences.  In our case, this hypothesis is automatically satisfied as $\Gamma_{(0)}(f,f)+\Gamma_{(1)}(f,f) = \left| \nabla f \right|^2$. For our purposes, we shall focus on their subRiemannian generalization of the Bonnet-Myers theorem.

\bgT[Baudoin-Garofalo]{BM}
If the generalized curvature inequality is satisfied with $\rho_1>0$  and the above hypothesis holds together with \rfE{com}  then $M$ is compact.
\enT

Combining this with \rfT{gci}, provides the following generalization of the examples considered in \cite{BaudoinGarofalo}

\bgT{BM2}
Under the same conditions as \rfT{gci}, if $\rho_1 >0$ then $M$ is compact.
\enT 
\section{Comparison with Riemannian curvatures}\setS{RM}

A common theme in the early development of subRiemannian geometry was the use of Riemannian approximations. More precisely, a Riemannian extension $g =g_0 \oplus g_1$ was chosen and then re-scaled as $g^\lambda = g_0 \oplus \lambda^2 g_1$. The behavior of these Riemannian metrics was then studied as $\lambda \to \infty$. The idea is that blowing up the vertical directions makes movement in these direction prohibitively expensive so the Riemannian geodesics should converge to the subRiemannian geodesics. Unfortunately, this is problematic for the study of the effects of curvature as this re-scaling makes the vertical curvatures much larger than the horizontal ones. However, useful information can be derived from this approach if instead we let $\lambda \to 0$.

In this section, we compute the Ricci and sectional curvatures of these scaled Riemannian metrics in terms of the basic connection. For simplicity, we shall specialize to the case where  $\dim VM=1$ and so the only basic grading applies.  We shall be able to provide alternative proofs to some of the results in \rfS{RC} and see the nature of the obstructions when the conditions are weakened.

To proceed, we fix a sRC-manifold $M$ and choose a Riemannian extension $g = g_0 \oplus g_1$. The basic connection will always be in terms of this metric. Throughout this section, $E_1,\dots,E_d$ will represent an orthonormal frame for $HM$ with respect to $g$ and $U$ will represent a unit length vector in $VM$, again with respect to $g$.

We refine the $J$ operator introduced earlier by defining
\bgE{J} \begin{split}  \aip{J^1(A,B)}{C}{} &= \aip{\tor(A,C)}{B_1}{},\\  \aip{J^0(A,B)}{C}{} &= \aip{ \tor(A,C)}{B_0}{}\end{split}\enE





\bgL{Comp}
For any sRC-manifold (with no restriction on $\dim VM$) the Levi-Civita connection associated to $g$ can be computed from the basic connection for $g$ as follows
\bgE{Comp}
\begin{split}
\bn_X Y &= \nabla_X Y  - \frac{1}{2} \tor(X,Y)  + J^1(X,Y)\\
\bn_T T  &= \nabla_T T  - \frac{1}{2} J^0(T,T)\\
\bn_T X &= \nabla_T X + \frac{1}{2}  J^0(X,T)  - \tor(T,X)_1\\
\bn_X T &= \nabla_X T  + \frac{1}{2} J^0(X,T)  -  \tor(X,T)_0
\end{split}
\enE
\enL

From these it is a straightforward, if brutal, computation to show that

\bgC{RmComp}
If $X,Y$ are horizontal vector fields and $T$ is a vertical vector field then
\bgE{XY}
\begin{split}
\bR(X,Y,Y,X) &= \srm(X,Y,Y,X)  - \frac{3}{4}  \left| \tor(X,Y) \right|^2  \\
& \qquad - \aip{J^1(Y,Y)}{J^1(X,X)}{} +\left| J^1(X,Y) \right|^2
\end{split}
\enE
\bgE{XT}
\begin{split}
\bR(T,X,X,T) &=  \srm(T,X,X,T) +\frac{1}{4}  \left|J^0(X,T)\right|^2  \\
& \qquad +  \aip{  \nabla \tor(T,X,X)  - \tor(X,\tor(X,T)) }{T}{} \\
& \qquad +\aip{\nabla \tor(X,T,T) }{X}{} - \left| \tor(X,T)_0 \right|^2
\end{split}
\enE
\bgE{XYT}
\begin{split}
\bR(X,Y,T,X) &=  \srm(X,Y,T,X) + \frac{1}{2} \aip{\nabla \tor (Y,X,X)}{T}{} \\
& \qquad   +\aip{ \nabla \tor( X,T,Y) - \nabla\tor (Y,T,X) }{X}{}
\end{split}
\enE
\enC
While this is far from a complete list of curvature terms, if we use properties of both Riemannian and subRiemannian curvatures and polarization identities, it is sufficient to compute all sectional and Ricci curvatures for the case $\dim VM=1$.

\bgR{neg}
If $M$ is strictly normal, then  \rfE{XT} reduces to
\[ \bR(T,X,X,T)=\frac{1}{4}  \left|J^0(X,T)\right|^2 \]
and so, if $HM$ bracket generates at step $2$,  there will always be at least one plane with postive sectional curvature. This means that a Riemannian approach to generalizing results concerning negative sectional curvatures is likely to be very difficult.
\enR


Provided that we only use constants for our re-scaling, it is easy to verify that the covariant derivatives for the basic connection associated to the re-scaled metric are unchanged from the base metric. Thus, paying careful attention to how each term scales, we can compute the Riemannian Ricci curvatures for the metrics $g^\lambda =g_0 \oplus \lambda^2 g_1$.

  For $Y \in HM$ and $T \in VM$, with inner products and norms computed in the unscaled metric
  \bgE{RcFla}
\begin{split}
\overline{\text{Rc}}^{\lambda} (Y,Y) &= \lambda^0 \left[ Rc(Y,Y) +  \aip{ \nabla \tor(U,Y,Y)}{U}{}    -\aip{\Tor(Y,Y,U)}{U}{} \right]  \\
& \qquad +\lambda^2 \left[  - \frac{1}{2} \sum\limits_i \left| \tor(E_i,Y) \right|^2 \right] \\
& \qquad + \lambda^{-2} \Big[  \aip{\nabla \tor(Y,U,U)}{Y}{}- \left| \tor(Y,U)_0 \right|^2 \\
&\qquad  + \sum\limits_i \left(  \left| J^1(E_i,Y)\right|^2 - \aip{J^1(E_i,E_i)}{J^1(Y,Y)}{} \right)  \Big] 
\end{split}
\enE
 \bgE{RcFlb}
\begin{split}
\overline{\text{Rc}}^{\lambda} (Y,T)  & =  \lambda^0 \left[  \sum\limits_i \aip{ \nabla \tor( E_i,T,Y) - \nabla\tor (Y,T,E_i) }{E_i}{}\right] \\
& \qquad + \lambda^2 \left[ \frac{1}{2} \aip{  \text{tr}_0 \nabla \tor (Y) }{T}{}  \right] \\
\end{split}
\enE
 \bgE{RcFlc}
\begin{split}
\overline{\text{Rc}}^{\lambda} (T,T) & = \lambda^0 \left[\sum _i \aip{\nabla \tor(E_i,T,T) }{E_i}{} - \left| \tor(E_i,T)_0 \right|^2 \right] \\
& \qquad + \lambda^2 \left[ \sum_i  \aip{  \nabla \tor(T,E_i,E_i)  -\Tor(E_i,E_i,T)}{T}{} \right] \\
& \qquad  + \lambda^4 \left[ \sum_i \frac{1}{4}  \left|J^0(E_i,T)\right|^2  \right] \\
\end{split}
\enE

For the case of a strictly normal sRC-manifold, these formulae greatly simplify to 
\bgE{RcFlSN}
\begin{split}
\overline{\text{Rc}}^{\lambda} (Y,Y) &= \lambda^0 Rc(Y,Y)  -\frac{\lambda^2}{2} \sum\limits_i \left| \tor(E_i,Y) \right|^2 \\
\overline{\text{Rc}}^{\lambda} (Y,T)  &=   \frac{\lambda^2}{2} \aip{  \text{tr}_0 \nabla \tor (Y) }{T}{} \\
\overline{\text{Rc}}^{\lambda} (T,T) & = \frac{\lambda^4}{4}  \sum_i  \left|J^0(E_i,T)\right|^2 = \frac{\lambda^4}{4} \sum\limits_{i,j} \left| \tor(E_i,E_j) \right|^2 
\end{split}
\enE
and so
\bgE{BGtensorRm}
 \begin{split}    \mathcal{R}(T+Y,T+Y) &=  \lim\limits_{\lambda \to 0} \overline{\text{Rc}}^{\lambda}   (Y+\lambda^{-2} T,Y+\lambda^{-2} T)  
  \end{split} 
 \enE

Next we note that if $T$ is unit length with respect to the base metric then for any smooth function
\[ \nabla^\lambda f = \nabla_0 f + \lambda^{-2} (Tf) T \]
which means that the Baudoin-Garofalo tensor applied to $\nabla f$ can expressed as a limit of Riemannian Ricci curvatures as follows:
\bgE{BGT}
 \mathcal{R}(\nabla f,\nabla f) =  \lim\limits_{\lambda \to 0} \overline{\text{Rc}}^{\lambda}   (\nabla^\lambda f, \nabla^\lambda f)  
\enE

\bgT{gci}
Under the same conditions as \rfT[RC]{gci} and the added assumption that $\dim VM=1$, there are constants $\lambda,c>0$ such that
\[  \overline{\text{Rc}}^{\lambda} (A,A) \geq c g^\lambda(A,A) \]
for all vectors $A$.
\enT

\pf Split $A = A_0+A_1$ and then note that
\begin{align*}
  \overline{\text{Rc}}^{\lambda} (A,A) &= \mathcal{R}(A_0 + \lambda^2 A_1,A_0+\lambda^2 A_1) -\frac{\lambda^2}{2} \sum\limits_i \left| \tor(E_i,A_0) \right|^2 \\
  & \geq \left(\rho_1 - \frac{\lambda^2 \kappa}{2} \right) g^\lambda(A_0,A_0) + \rho_2 \lambda^2 g^\lambda(A_1,A_1) 
  \end{align*}
  Thus for very small $\lambda>0$, we can take $c= \min\left\{ \lambda^2\rho_2, \rho_1 - \frac{\lambda^2 \kappa}{2} \right\} >0$.
  
\epf

Combining this with the classical Myers theorem yields

\bgC{gci}
Under the same conditions as \rfT[RC]{gci} and the added assumption that $\dim VM=1$, $M$ is compact and has finite fundamental group.
\enC

\bgR{BG}
If $\dim VM>1$ this characterization of the Baudoin-Garofalo tensor will fail as there will be additional vertical Ricci curvature terms. One real advantage of the heat kernel method established in \cite{BaudoinGarofalo} is that the result still works in this case. Rather than working with the entire the Riemannian  Ricci curvature, the Baudoin-Garofalo tensor focuses on just the portion that is actually needed to ensure compactness. The result is thus stronger than that which could be obtained using purely Riemannian methods. Indeed if $\dim VM=2$, then positivity of the Baudoin-Garofalo tensor will then impose topological constraints on $VM$. As $VM$ is integrable it will generate a foliation. If the conditions of the previous theorem are met then the leaves of this foliation cannot be non-compact embedded submanifolds. Thus either they must be immersed or their Gaussian curvature must be non-negative.
\enR

\bgR{est}
It should be noted, that \rfC{gci} does not come with a diameter estimate. The classical Bonnet-Myers theorems yield diameter estimates in terms of the distance associated to the metrics $g^\lambda$. However for small $\lambda$, this distance has little relation to the subRiemannian distance. The heat kernel methods underlying \rfT[RC]{gci} do produce diameter estimates in terms of "\textit{carr\'e du champ}" distance (see \cite{BaudoinGarofalo}) which is much more closely related to the subRiemannian distance.
\enR

If we do not restrict to the strictly normal case, then this Riemannian approach immediately has problems. If we send $\lambda \to \infty$, we see that $\overline{Rc}^\lambda(Y,Y) \to -\infty$, so any Riemannian results for positive Ricci curvature will immediately be lost. Since there are very few topological consequences of negative Ricci curvature, this approach is unlikely to bear fruit. If however we let $\lambda \to 0$, then we run into the issue that the subRiemannian Ricci curvature for the horizontal terms isn't the dominant term. Instead we must deal with the symmetric 2-tensors
\bgE{B}
\begin{split}
\mathcal{B}(X,Y) &=  \aip{\nabla \tor(X,U,U)}{Y}{}- \aip{ \tor(X,U)_0}{\tor(Y,U)_0 }{} \\
\mathcal{K}(X,Y) &=  \sum\limits_i \left(  \aip{J^1(E_i,X)}{ J^1(E_i,Y)}{}  - \aip{ J^1(E_i,E_i)}{J^1(X,Y)}{} \right)
\end{split}
\enE
where again $U$ is a unit length vertical vector. The tensor $\mathcal{B}$ is a genuine sRC-invariant when $\dim VM=1$, but has no good invariant generalization when $\dim VM>1$.  
However $\mathcal{K}$ is  only a vertically conformal sRC-invariant.

With these caveats in mind, we do however obtain the following theorem

\bgT{BadBonnet}
Let $M$ be an sRC-manifold with $\dim VM=1$ and bounded curvature and torsion. If there are constants $a,b > 0$ such that for all horizontal vectors $Y$,
\bgE{Bad}
\begin{split}
 \text{tr}_0 \mathcal{B} &\geq a\\
\mathcal{B}(Y,Y) +\mathcal{K}(Y,Y)  &\geq b \left|Y\right|^2 
\end{split}
\enE
then $M$ is compact and has finite fundamental group.
\enT

\pf
The condition of bounded curvature implies that for small $\lambda$ there will be some, possibly large,  constant $M$ such that
\[  2\overline{Rc}^\lambda (T,Y) \leq 2M \left| T \right| \left| Y \right| \leq  \frac{a}{4} \left| T\right|^2 + \frac{4M^2}{a} \left|Y\right|^2. \]
Since $ \text{tr}_0 \mathcal{B} \geq a$ globally, for sufficiently small $\lambda$, we will have
\[  \overline{Rc}^\lambda (T,T)  \geq \frac{a}{2} \left|T\right|^2 \]
and since  $\mathcal{K}(Y,Y) +\mathcal{B}(Y,Y)  \geq b \left|Y\right|^2 $, again for small $\lambda$, we have
\[  \overline{Rc}^\lambda (Y,Y) \geq \frac{b}{2 \lambda^2}  \left|Y\right|^2 \]
But then for small enough $\lambda$
\[ \overline{Rc}^\lambda (T+Y,T+Y) \geq  \left(  \frac{b}{2\lambda^2} - \frac{4M^2}{a} \right) \left|Y\right|^2 + \frac{a}{4} \left|T\right|^2 \]
For very small $\lambda$, both coefficients will be postive, so 
\[  \overline{Rc}^\lambda (T+Y,T+Y) \geq c \left|T+Y\right|^2 \]
for some positive constant $c$. The result then follows from the classical Myers theorem.

\epf

This is a purely subRiemmanian result as the conditions are trivially false when restricted to Riemannian manifolds. However, it is somewhat  unsatisfactory in nature. It would seem reasonable to conjecture that for sRC-manifolds (or at least those that are in some sense nearly strictly normal ) that there would be some sort of analogue of \rfT[RC]{BM2} where the dominant terms are genuine subRiemannian Ricci tensors. However, it appears that  to prove it will be necessary to create new subRiemannian techniques such as the heat kernel methods of \cite{BaudoinGarofalo} rather than fall back on existing Riemannian methods. The author expects the basic connection developed to provide a solid computational foundation for such techniques.

\bibliographystyle{plain}
\bibliography{References}

\end{document}